\documentclass[12pt]{article}
\usepackage{amssymb}
\newcommand{\ds}{\displaystyle}
\renewcommand{\a}{{\alpha}}
\renewcommand{\b}{{\beta}}
\renewcommand{\d}{{\delta}}
\newcommand{\g}{{\gamma}}
\newcommand{\tfrac}[2]{{\mbox{$\frac{#1}{#2}$}}}
\newcommand{\proofend}{\hfill{$\Box$}}
\newcommand{\bbN}{{\mathbb N}}

\newcommand{\bde}{{\mathbf{e}}}
\newcommand{\bdx}{{\mathbf{x}}}
\newcommand{\bdy}{{\mathbf{y}}}

\newtheorem{satz}{Satz}[section]
\newtheorem{theorem}{Theorem}[section]
\newtheorem{pro}[satz]{Proposition}
\newtheorem{cor}[satz]{Corollary}
\newtheorem{lem}[satz]{Lemma}

\begin{document}
\title{Ces\`aro means of Jacobi expansions on the parabolic biangle}
\author{W. zu Castell\thanks{Institute of Biomathematics and
Biometry, Helmholtz Zentrum M\"{u}nchen, German Research Center for
Environmantal Health, Ingolst\"{a}dter  Landstra\ss e 1,
85764 Neuherberg, Germany. [castell,filbir]@helmholtz-muenchen.de},
F. Filbir{$^*$},
Y. Xu\thanks{Department of Mathematics, University of Oregon, Eugene, 
Oregon 97403-1222, U.S.A.\newline yuan@uoregon.edu . This work was
partially supported by the NSF grant DMS-0604056.} }
\date{\today}

\maketitle

{\small\bf Abstract. }
{\small 
We study Ces\`aro $(C,\delta)$ means for two-variable Jacobi polynomials
on the parabolic biangle 
$B=\{(x_1,x_2)\in{\mathbb R}^2:0\leq x_1^2\leq x_2\leq 1\}$. Using the
product formula derived by Koornwinder \& Schwartz for this polynomial
system, the Ces\`aro operator can be interpreted as a convolution operator.
We then show that the Ces\`aro $(C,\delta)$ means of the orthogonal
expansion on the biangle are uniformly bounded if $\delta>\alpha+\beta+1$,
$\alpha-\frac 12\geq\beta\geq 0$. Furthermore, for 
$\delta\geq\alpha+2\beta+\frac 32$ the means define positive linear
operators.
}

{\small\bf Keywords. }
{\small orthogonal expansion, Ces\`aro summability, parabolic biangle, 
two-variable orthogonal polynomials, positive linear operators, 
convolution operators}

\section{Introduction}
\label{sec:introduction}
\setcounter{equation}{0}

Multivariate analogs of classical orthogonal polynomials are
of great interest in many areas of applied analysis and approximation.
Although, there is a beautiful theory of orthogonal polynomials in
one variable (cf. \cite{Sz75}), it is much harder to develop a
picture as comprehensive as in the univariate case for orthogonal
polynomials in several variables (cf. \cite{Du01}). Among classical
orthogonal polynomials in one variable the family of Jacobi polynomials 
plays a special role. Their importance is partly due to the fact that 
there are connections of these polynomials to group representations and
the eigensystem of the Laplace-Beltrami operator for certain compact
symmetric spaces. In the two-dimensional setting, Koornwinder
\cite{Ko75,Ko75a} discusses seven classes of Jacobi polynomials in two
variables, derived by either expressions in terms of univariate 
Jacobi polynomials, studying spherical functions on homogeneous
spaces of rank 2, or using quadratic transformations of known
examples of Jacobi polynomials, as well as analytic continuation
with respect to some parameters. Among these classes, there are
bivariate Jacobi polynomials on the biangle. The \emph{parabolic biangle} 
is the closed subset of $\mathbb{R}^2$ defined by
$$
B=\{\bdx=(x_1,x_2)\in\mathbb{R}^2: 0\leq x_1^2\leq x_2\leq 1\}.
$$
The system of bivariate orthogonal polynomials defined on the set $B$
has first been introduced by Agahanov \cite{Ag65}. The polynomials can be 
explicitly expressed in terms of the univariate Jacobi polynomials. In 
order to provide the explicit formulae let us recall some basic facts for 
Jacobi polynomials. 

\medskip
For $\alpha,\beta>-1$ and $n\in\mathbb{N}_0$ the Jacobi polynomials are 
defined as 
$$
P_n^{(\alpha,\beta)}(x)
\, = \,
\frac{(\alpha+1)_n}{n!}\, 
F\left(-n,n+\alpha+\beta+1;\alpha+1;\frac{1-x}{2}\right),
$$
where $F(a,b;c;x)$ denotes the Gaussian hypergeometric function. 
Given the weight function $w^{(\alpha,\beta)}$ defined as
$$
w^{(\alpha,\beta)}(x)
=\frac{\Gamma(\alpha+\beta+2)}
      {2^{\alpha+\beta+1}\Gamma(\alpha+1)\Gamma(\beta+1)}
(1-x)^\alpha(1+x)^\beta,
\quad x\in [-1,1],
$$
the Jacobi polynomials satisfy the following orthogonality
relation.
$$
\int_{-1}^1 P_n^{(\alpha,\beta)}(x)P_m^{(\alpha,\beta)}(x)
w^{(\alpha,\beta)}(x)\, dx
\, = \,
\delta_{n,m}h_n^{(\alpha,\beta)},
$$
where
\begin{equation}
\label{haar-weights}
h^{(\alpha,\beta)}_n
\, = \, 
\frac{(\alpha+1)_n(\beta+1)_n(\alpha+\beta+1)}
     {n!(\alpha+\beta+2)_n(\alpha+\beta+2n+1)_n}.
\end{equation}
We are now ready to define the orthogonal polynomials system on $B$. 
For $\alpha,\beta>-1$ and $n\geq k\geq 0$ let 
\begin{equation}
\label{def-orthpol}
P^{\alpha,\beta}_{n,k}(\bdx)
\, = \,
P_{n-k}^{(\alpha-\frac 12,\beta+k)}(2x_2-1)\cdot x_2^{\frac{k}{2}}\, 
 P_{k}^{(\beta-\frac 12,\beta-\frac 12)}(x_1x_2^{-\frac{1}{2}}),
\quad \bdx=(x_1,x_2)\in B.
\end{equation}
Note that the polynomial is of degree $k$ in $x_1$ and of degree $n-k$ 
in $x_2$ and thus, of total degree $n$. Furthermore, 
$$
P^{\alpha,\beta}_{n,k}(\bde)
\, = \,
\frac{\left(\alpha+\frac 12\right)_{n-k}\left(\beta+\frac 12\right)_k}
     {(n-k)!k!},
$$
where $\bde=(1,1)\in B$. Moreover, for $\alpha,\beta$ fixed these 
polynomials are orthogonal with respect to the weight function 
\begin{equation}
\label{eq-weightW}
W^{\alpha,\beta}(\bdx)
\, = \,
\frac{\Gamma\left(\alpha+\beta+\frac 32\right)}
     {\Gamma\left(\frac 12\right)
      \Gamma\left(\alpha+\frac 12\right)\Gamma\left(\beta+\frac 12\right)}\,
(1-x_2)^{\alpha-\frac 12}(x_2-x_1^2)^{\beta-\frac 12},
\quad \bdx=(x_1,x_2)\in B.
\end{equation}
Thus, for the $L^2-$norm of the polynomials $P_{n,k}^{\alpha,\beta}$ we 
obtain the expression
\begin{equation}
\label{eq-2norm}
\begin{array}{ll}
\big(g^{\alpha,\beta}_{n,k}\big)^{-1}
& = \ds\int_B |P_{n,k}^{\alpha,\beta}(\bdx)|^2\, 
W^{\alpha,\beta}(\bdx)\ d\bdx\\[3ex]
& = \ds
h_{n-k}^{(\alpha-\frac 12,\beta+k)}h_k^{(\beta-\frac 12,\beta-\frac 12)}.
\end{array}
\end{equation}
To keep the notation simple we write $L^p(B)$ for the Lebesgue spaces 
$L^p(B,W^{\alpha,\beta})$, $1\leq p\leq\infty$, where for $p=\infty$ 
the space $L^\infty(B)$ is understood to be the space of continuous 
functions on $B$ with the supremum norm. 

We want to study Ces\`aro means of orthogonal expansions with respect
to the system $\{P_{n,k}^{\alpha,\beta}: 0\leq k\leq n\}$, $n\in\bbN_0$,
on $B$.  

\medskip
For $\delta\geq 0$, the \emph{Ces\`aro $ (C,\delta)$ means} $s_n^\delta$
of a sequence $(c_n)_{n=0}^\infty$ are defined as
\begin{equation}
\label{cesaro_means}
s_n^\delta
\, = \,
\frac 1{{{n+\delta}\choose{n}}}\,
\sum_{k=0}^n {{n-k+\delta-1}\choose{n-k}}\, s_k
\, = \,
\frac 1{{{n+\delta}\choose{n}}}\,
\sum_{k=0}^n {{n-k+\delta}\choose{n-k}}\, c_k,
\end{equation}
where $s_n$ denotes the $n$-th partial sum $\sum_{k=0}^n c_k$. The sequence 
$(s_n)_{n=0}^\infty$ is called \emph{$(C,\delta)$ summable} if $s_n^\delta$ 
converges for $n\to\infty$. For simplicity, we set 
$A_n^\delta={{n+\delta}\choose n}$.

\medskip
The study of Ces\`aro means for orthogonal expansions has a long history.
It is worthwhile to note that in the multivariate setting, there is a
strong influence of the domain of definition, which is up to now not fully
understood. For example, the Ces\`aro means for the 
\emph{Chebyshev weight function} $1/\sqrt{(1-x^2)(1-y^2)}$ converge
uniformly on the unit square for all $\delta>0$. On the unit ball, there is 
in contrast a \emph{critical index} for the weight $1/\sqrt{1-\|x\|^2}$, which 
means, uniform convergence does not hold if $\delta\leq\frac{d-1}2$. The 
same is true for the standard simplex (cf. \cite{DaiX, Li03} for further details).
The parabolic biangle is in some sense an intermediate region having both, 
a smooth curved boundary as well as singularities. As we will show below, 
for $1/ \sqrt{(1-x_2)(x_2-x_1^2)}$ on this region uniform convergence holds 
if $\delta>1$.

\medskip
In the following section we will introduce the Ces\`aro means for the
orthogonal expansion on the biangle. Since one crucial point for our
following discussion is the fact that the Ces\`aro means can be
interpreted as a convolution operator, we recall some facts on the
convolution structure on the parabolic biangle.
In Section \ref{sec:summability} our main results are stated, while
the proof of Theorem \ref{th-main}, our main result, is given in the 
last section. Throughout the paper we will use the letter $c$ for a 
generic constant which does not have to be the same in every occurrence.

\section{Ces\`aro means on the biangle} 
\label{sec:orthogonal_expansions}
\setcounter{equation}{0}

The Fourier coefficient of a function $f\in L^1(B)$ with respect to the
orthogonal system $\{P_{n,k}^{\alpha,\beta}:0\leq k\leq n\}$, 
$n\in\bbN_0$, is given by 
\begin{equation}
\widehat{f}(n,k)=
\int_B f(\mathbf{y})\, P_{n,k}^{\alpha,\beta}(\mathbf{y})\ 
W^{\alpha,\beta}(\mathbf{y})\ d\mathbf{y},
\end{equation}
where $0\leq k\leq n$, $n\in\bbN_0$.

The projection of $f$ onto the space of polynomials of degree $n$ is 
given by 
\begin{equation}
\label{eq-partial}
\mathcal{P}_nf(\mathbf{x})
\, =\,
\sum_{k=0}^n \widehat{f}(n,k)\, P_{n,k}^{\alpha,\beta}(\mathbf{x})
 g^{\alpha,\beta}_{n,k}
\, = \,
\int_B f(\mathbf{x})
P_n(\mathbf{x},\mathbf{y})\, W^{\alpha,\beta}(\mathbf{y})\ d\mathbf{y},
\end{equation}
with kernel 
$$
P_n(\mathbf{x},\mathbf{y})=
\sum_{k=0}^nP_{n,k}^{\alpha,\beta}(\mathbf{x})P_{n,k}^{\alpha,\beta}
(\mathbf{y})g^{\alpha,\beta}_{n,k}.
$$
We want to study approximate expansions of functions using summability 
methods. Let us therefore recall the definition of Ces\`aro means.

The Ces\`aro $(C,\delta)$ means of the expansion (\ref{eq-partial}) are 
given by 
\begin{equation}
\label{eq-cesaro3}
\mathcal{S}^\delta_nf(\mathbf{x})
\, = \,
\sum_{k=0}^n \frac{A_{n-k}^\delta}{A_n^\delta}\,\mathcal{P}_k f(\mathbf{x})
\, =\,
\int_B f(\mathbf{y})
\mathcal{K}^\delta_n(\mathbf{x},\mathbf{y})\, 
W^{\alpha,\beta}(\mathbf{y})\ d\mathbf{y},
\end{equation}
where the summability kernel is given by
\begin{equation}
\label{eq-kernel1}
\begin{array}{ll}
\mathcal{K}^\delta_n(\mathbf{x},\mathbf{y})
& =
\ds\sum_{k=0}^n \frac{A_{n-k}^\delta}{A_n^\delta}\, P_k(\mathbf{x},\mathbf{y})
\, = \,
\sum_{k=0}^n\sum_{l=0}^k 
\frac{A_{n-k}^\delta}{A_n^\delta}\,
P_{n,k}^{\alpha,\beta}(\mathbf{x})P_{n,k}^{\alpha,\beta}(\mathbf{y}) 
g^{\alpha,\beta}_{n,k}\\[3ex]
& = 
\ds\sum_{k=0}^n\frac{A_{n-k}^{\delta-1}}{A_n^\delta}\,
\mathcal{K}_k(\mathbf{x},\mathbf{y}).
\end{array}
\end{equation}
Here and for the remaining part of the paper we simply write $\mathcal{S}_n$ for 
$\mathcal{S}_n^0$ and $\mathcal{K}_n$ for $\mathcal{K}_n^0$. 
Now, a key observation is that the operators $\mathcal{P}_n$ and 
$\mathcal{S}_n^\delta$ can be written as convolution operators. In order to 
see how this works, we need to outline some facts concerning the convolution 
structure on $B$ associated with the polynomial 
system $\{P_{n,k}^{\alpha,\beta}: 0\leq k\leq n\}$. 

Convolution structures are commonly based on product formulas. For the 
polynomials $P_{n,k}^{\alpha,\beta}$ such a construction was given by Koornwinder 
and Schwartz \cite{Ko97}. 
We need some notation. 
For $\bdx=(x_1,x_2)$ with $|x_1|, |x_2|\in [0,1]$, $0\leq r\leq 1$, and 
$\psi,\psi_j\in [0,\pi], j=1,2,3$, define
$$
\begin{array}{rl}
D(\bdx;r,\psi)=&
x_1x_2+(1-x_1^2)^{1/2}(1-x_2^2)^{1/2}\, r\cos\psi ,\\[1ex]
E(\bdx;r,\psi)=&
(x_1^2x_2^2+(1-x_1^2)(1-x_2^2)r^2 + 
2x_1x_2(1-x_1^2)^{1/2}(1-x_2^2)^{1/2}\, r\cos\psi)^{1/2},
\end{array}
$$
and, setting $\bdy=(y_1,y_2)$,
\begin{eqnarray*}
F(\bdx,\bdy;r,\psi_1,\psi_2,\psi_3) 
 =
E(\bdx;r,\psi_1)
D\left[\frac{D(\bdx;r,\psi_1)}{E(\bdx;r,\psi_1)},
       D\left(\frac{x_1}{x_2},\frac{y_1}{y_2};1,\psi_2\right);1,\psi_3
\right],
\end{eqnarray*}
where $\bdx,\bdy\neq 0$. 
Moreover, we define the following measures 
$$
dm^{\alpha,\beta}(r,\psi)
\, = \,
\frac{2\Gamma(\alpha+1/2)\Gamma(\beta+1/2)}
     {\Gamma(\alpha-\beta)\Gamma(\beta+1/2)\Gamma(\beta)\Gamma(1/2)}\,
     (1-r^2)^{\alpha-\beta-1}
      r^{2\beta}\sin^{2\beta-1}\psi\ dr\, d\psi
$$
and
\begin{eqnarray*}
d\mu^{\alpha,\beta}(r,\psi_1,\psi_2,\psi_3)
& = &
c_{\alpha,\beta}(1-r^2)^{\alpha-\beta-3/2}r^{2\beta} \\
& & \times \,
\sin^{2\beta-1}\psi_2\ 
\sin^{2\beta-1}\psi_3\ \sin^{2\beta}\psi_1\ drd\psi_1\, d\psi_2\, d\psi_3,
\end{eqnarray*}
where $c_{\alpha,\beta}$ is a constant so that $d\mu^{\alpha,\beta}$ is a 
probability measure. 

Koornwinder and Schwartz \cite{Ko97} proved the following product 
formula for the orthogonal polynomials on the parabolic biangle. 

\begin{theorem}
Let $\alpha-\frac 12\geq\beta\geq 0$. Assume that $0\leq|x_1|\leq x_2\leq 1$ and 
$0\leq|y_1|\leq y_2\leq 1$. Then if 
$\bdx=(x_1,x_2),\bdy=(y_1,y_2)\in B\setminus\{0\}$, we 
have that
\begin{eqnarray*}
\lefteqn{
P_{n,k}^{\alpha,\beta}(x_1,x_2^2)\, P_{n,k}^{\alpha,\beta}(y_1,y_2^2) \, = }\\
& &
P_{n,k}^{\alpha,\beta}(\bde)\,
\int_{[0,1]\times [0,\pi]^3}P_{n,k}^{\alpha,\beta}
(E^2(\bdx;r,\psi),F(\bdx,\bdy;r,\psi_1,\psi_2,\psi_3))\ 
d\mu^{\alpha,\beta}(r,\psi_1,\psi_2,\psi_3).
\end{eqnarray*}
If $0\leq|x_1|\leq x_2\leq 1$ and $\bdy=0$ we have that
$$
P_{n,k}^{\alpha,\beta}(x_1,x_2^2)\, P_{n,k}^{\alpha,\beta}(0)=
P_{n,k}^{\alpha,\beta}(\bde)\,
\int_{[0,1]\times [0,\pi]^3}P_{n,k}^{\alpha,\beta}
(E^2(\bdx;r,\psi),D(\bdx;r,\psi_1))\ 
dm^{\alpha,\beta+1/2}(r,\psi_1).
$$
\end{theorem}

For our aim we do not need the product formula in its explicit form, rather than 
its existence, which gives rise to the convolution structure on the parabolic 
biangle. Thus, let us restate the product formula in a convenient form: 
%
\begin{equation}
\label{PF-1}
P_{n,k}^{\alpha,\beta}(\mathbf{x})\, P_{n,k}^{\alpha,\beta}(\mathbf{y})
\, = \,
P_{n,k}^{\alpha,\beta}(\bde)\,
\int_{B}P_{n,k}^{\alpha,\beta}(\mathbf{z})\ 
d\omega_{\mathbf{x},\mathbf{y}}(\mathbf{z}),
\end{equation}
where $d\omega_{\mathbf{x},\mathbf{y}}(\mathbf{z})$ is a probability measure.

Formula (\ref{PF-1}) gives rise to a generalized translation operator 
on $B$ in the following way. For $f\in C(B)$ we define 
$$
T_{\mathbf{x}}f(\mathbf{y})
\, = \,
\int_{B}f(\mathbf{z})\ d\omega_{\mathbf{x},\mathbf{y}}(\mathbf{z}),
\quad \bdx,\bdy\in B.
$$
It can be shown that $T_{\mathbf{x}}$ extends to a bounded linear operator 
on $L^p(B)$ with $\|T_\bdx\|_p\leq 1$, $\bdx\in B$, 
for every $1\leq p\leq\infty$. 
The convolution of functions $f,g\in L^1(B)$ is then defined as 
\begin{equation}
\label{conv}
f\ast g(\mathbf{x})
\, = \,
\int_B f(\mathbf{y})T_{\mathbf{x}}g(\mathbf{y})\ 
W^{\alpha,\beta}(\mathbf{y})\, d\mathbf{y},
\quad \bdx\in B.
\end{equation} 
This convolution product is associative and commutative. Moreover, the following 
estimate holds true for all $f,g\in L^p(B)$: 
\begin{equation}\label{norm}
\|f\ast g\|_p\leq\|f\|_1\, \|g\|_p,\quad 1\leq p\leq\infty.
\end{equation}
In view of formulae (\ref{PF-1}) and (\ref{conv}) it becomes  obvious that the 
operator $\mathcal{S}_n^\delta$ can be written as convolution operator. To be
precise, we have the following result.

\begin{pro}
If $\alpha-\frac 12\geq\beta\geq 0$, $\delta >0$, and $\bdx\in B$, then for
every $f\in L^p(B)$ we have that
\begin{equation}
\label{cesaro_as_convolution}
\mathcal{S}_n^\delta f(\mathbf{x}) \, = \, K_n^\delta\ast f(\mathbf{x}),
\end{equation}
where $K_n^\delta(\mathbf{x})= \mathcal{K}_{n}^\delta(\mathbf{x},\mathbf{e})$
and $\mathbf{e}=(1,1)$. 
\end{pro}

\medskip
{\it Proof. }
From equations (\ref{eq-cesaro3}) and (\ref{eq-kernel1}) we obtain, using the
product formula (\ref{PF-1}), that
\begin{eqnarray*}
\mathcal{S}_n^\delta f(\bdx) 
& = &
\int_B f(\bdy)\,
\sum_{k=0}^n\sum_{l=0}^k \frac{A^\delta_{n-k}}{A_n^\delta}\,
P_{n,k}^{\alpha,\beta}(\bde) T_\bdx P_{n,k}^{\alpha,\beta}(\bdy)\,
W^{\alpha,\beta}(\bdy)\, d\bdy \\
& = & 
K_n^\delta \ast f(\bdx).
\end{eqnarray*}
\proofend

Now from inequality (\ref{norm}) an estimate for the operator norm of  
$\mathcal{S}_n^\delta$ follows, i.e., 
$\Vert\mathcal{S}_n^\delta\Vert_p\leq \|K_n^\delta\|_1,\ 1\leq p\leq\infty$.
Notice that $K_n^\delta(\bdx)$ is the kernel $ \mathcal{K}_{n}^\delta(\mathbf{x},
\mathbf{y})$ when $\mathbf{y}= \mathbf{e}$. The uniform convergence of 
the $(C,\delta)$ means on $B$ is reduced to convergence at a point. We state
this formerly as a corollary.

\begin{cor} 
\label{cor2.2}
If $\alpha-\frac 12\geq\beta\geq 0$, $\delta>0$, and $f\in L^p(B)$,
the Ces\`aro $(C,\delta)$ means $\mathcal{S}_n^\delta f$ converge uniformly
on $B$ if $\mathcal{S}_n^\delta f$ converges at the point $\mathbf{e}$, which 
holds, in turn, if 
\begin{equation} \label{L1norm}
\|K_n^\delta\|_1 = \int_B |K_n^\delta(\bdx)| W^{(\alpha,\beta)}(\bdx) d\bdx \le c < \infty
\end{equation}
uniformly in $n$.
\end{cor} 

\section{Summability of Ces\`aro expansions on the biangle} 
\label{sec:summability}
\setcounter{equation}{0}

As shown in Corollary \ref{cor2.2}, it is sufficient to establish (\ref{L1norm}). For
this purpose, it is essential to derive a closed formula for the kernel $K_n^\delta$,
which we state below. 

\begin{theorem}
\label{th-cfkernel}
For $\alpha,\beta>-\frac 12$ and $\bdx\in B$ we have that 
\begin{eqnarray}
\label{eq-cfkernel}
K_n(\mathbf{x})
& = &
P_n^{(\alpha+\beta+\frac{1}{2},\beta)}(1)\\
\nonumber
& &
\times\, 
\int_{-1}^1 P_n^{(\alpha+\beta+\frac{1}{2},\beta)}
(\frac{1}{2}(1+t)^2+(1-t^2)x_1+\frac{1}{2}(1-t)^2x_2-1)\ 
w^{(\alpha+\beta+\frac{1}{2},\beta)}(t)\ dt.
\end{eqnarray}
\end{theorem}

{\it Proof. } 
We derive this formula from the addition formula for Jacobi polynomials  
$P_n^{(\alpha,\beta)}$ established by Koornwinder (cf. \cite[(4.14)]{Ko75}).
\begin{equation}
\label{eq-add}
\begin{array}{l}
P_n^{(\alpha,\beta)}(\frac{1}{2}(1+\xi)(1+\eta)+\frac{1}{2}
(1-\xi)(1-\eta)r^2+(1-\xi^2)^{1/2}(1-\eta^2)^{1/2}\, r\cos\psi -1)\\[2ex]
\begin{array}{rl}
=\ds\sum_{k=0}^n\sum_{l=0}^ka^{(\alpha,\beta)}_{n,k,l} &
(1-\xi)^{(k+l)/2}(1+\xi)^{(k-l)/2}P_{n-k}^{(\alpha+k+l,\beta+k-l)}(\xi)\\[2ex]
\cdot&(1-\eta)^{(k+l)/2}(1+\eta)^{(k-l)/2}
P_{n-k}^{(\alpha+k+l,\beta+k-l)}(\eta)\\[2ex]
\cdot& P_{l}^{(\alpha-\beta-1,\beta+k-l)}(2r^2-1)\, r^{k-l}\, 
P_{k-l}^{(\beta-1/2,\beta-1/2)}(\cos\psi),
\end{array}
\end{array}
\end{equation}
where the coefficients $a^{(\alpha,\beta)}_{n,k,l} $ are given by
\begin{eqnarray}
\label{eq-koeffAF}
a_{n,k,l}^{(\alpha,\beta)}
& = &
(k+l+\alpha)(k-l+\beta) \\
\nonumber
& & \times \,
\frac{(n+\alpha+\beta+1)_k(2\beta+1)_{k-l}(n-l+\beta+1)_l
      (n-k)!}
     {2^{2k}(k+\alpha)(\frac{1}{2}(k-l)+\beta)(\beta+1)_k
      (k+\alpha+1)_{n-k+l}(\beta+1/2)_{k-l}}.
\end{eqnarray}
By definition (\ref{eq-kernel1}) we have that
$$
K_n(\bdx,\bde) 
\, = \, 
\sum_{k=0}^n \sum_{l=0}^k 
P_{k,l}^{\alpha,\beta}(\bdx)P_{k,l}^{\alpha,\beta}(\bde)
g_{n,k}^{\alpha,\beta}.   
$$
From the fact that $P_n^{(\alpha,\beta)}(1) = (\alpha+1)_n /n!$, 
it follows that 
$$
P_{k,l}^{\alpha,\beta}(\bde) 
\, = \,
P_{k-l}^{(\alpha-\frac12,\beta+k)}(1)   
  P_{l}^{(\beta-\frac12,\beta-\frac12)}(1) =  
 \frac{(\alpha+\frac12)_{k-l}(\beta+\frac12)_l}{(k-l)! l!}. 
$$    
In the addition formula (\ref{eq-add}), 
set $\xi = \eta =t$, $r = \sqrt{x_2}$, $r\cos\psi=x_1$, and replace 
$\alpha$ by $\alpha+\beta+1/2$ to obtain
\begin{eqnarray*}
\lefteqn{   
P_n^{(\alpha+\beta+\frac12,\beta)}\Big(\tfrac{1}{2} (1+t)^2 +   
  (1-t^2) x_1+ \tfrac{1}{2} (1-t)^2 x_2  -1\Big)
} \\ 
& = \displaystyle{\sum_{k=0}^n \sum_{l=0}^k} 
    a_{n,k,l}^{(\alpha+\beta+\frac12,\beta)} 
    (1-t)^{k+l} (1+t)^{k-l} 
    \Big[P_{n-k}^{(\alpha+\beta+k+l+\frac12,\beta+k-l)} (t)\Big]^2 \\
& \times \,
  P_l^{(\alpha-\frac12,\beta+k-l)}(2x_2-1)x_2^{\frac{k-l}{2}}  
  P_{k-l}^{(\beta-\frac12,\beta-\frac12)}\Big(\frac{x_1}{\sqrt{x_2}}\Big). 
\end{eqnarray*} 
The product of the last three terms in the sum is precisely 
$P_{k,k-l}^{\alpha,\beta}$. Integrating the above equation with respect to 
$w^{(\alpha+\beta+\frac12,\beta)}(t)dt$ gives 
\begin{eqnarray*}
& & \int_{-1}^1 P_n^{(\alpha+\beta+\frac12,\beta)}
\Big(\tfrac{1}{2} (1+t)^2 +  
  (1-t^2) x_1 +\tfrac{1}{2} (1-t)^2 x_2 -1
\Big) 
w^{(\alpha+\beta+\frac12,\beta)}(t)dt
\\
& & \, = \, \sum_{k=0}^n \sum_{l=0}^k  a_{n,k,l}^{(\alpha+\beta+\frac12,\beta)}
        b_{n,k,l}^{(\alpha,\beta)}  P_{k,k-l}^{\alpha,\beta}(x),\\
\end{eqnarray*}
where 
$$
b_{n,k,l}^{(\alpha,\beta)}
\, = \,
\int_{-1}^1 \Big[P_{n-k}^{(\alpha+\beta+k+l+\frac12,\beta+k-l)}(t)\Big]^2 
(1-t)^{k+l} (1+t)^{k-l}w^{(\alpha+\beta+\frac12,\beta)}(t)dt.  
$$
Using (\ref{haar-weights}) and the explicit formula of $a_{n,k,l}^{(\alpha+\beta+
\frac12,\beta)}$ it can be verified that 
$$
a_{n,k,l}^{(\alpha+\beta+\frac12,\beta)} b_{n,k,l}^{(\alpha,\beta)} =
g_{k-l,k}^{\alpha,\beta} P_{k,k-l}^{\alpha,\beta}(\bde) 
  [P_n^{(\alpha+\beta+\frac12,\beta)}(1)]^{-1}
  h_n^{(\alpha+\beta+\frac12,\beta)}.
$$
This proves the stated formula.
\proofend

We will use the abbreviation $z(\bdx;t)$ for the argument of the univariate 
Jacobi polynomial in equation  (\ref{eq-cfkernel}), i.e.,
$$
z(\bdx;t) 
\, = \,
\mbox{$\frac 12$}(1+t)^2+(1-t^2)x_1+\mbox{$\frac 12$}(1-t)^2x_2-1,
\qquad \bdx=(x_1,x_2)\in B. 
$$
Using Theorem \ref{th-cfkernel} we obtain the following form of the 
Ces\`aro kernel.

\medskip
\begin{cor}
\label{kernel-reduced}
Let $\alpha-\frac 12\geq\beta\geq 0$ and $\delta >0$. Then for all
$\bdx\in B$ and $f\in L^p(B)$, $1\leq p\leq\infty$,
we have that 
\begin{eqnarray}
\nonumber
K_n^\delta(\mathbf{x})
& = & 
\sum_{k=0}^n\frac{A_{n-k}^{\delta-1}}{A_n^\delta}\,
h_k^{(\alpha+\beta+\frac 12,\beta)}
P_k^{(\alpha+\beta+\frac 12,\beta)}(1)\,
\int_{-1}^1 P_k^{(\alpha+\beta+\frac{1}{2},\beta)}(z(\bdx;t))\, 
w^{(\alpha+\beta+\frac{1}{2},\beta)}(t)\, dt\\
\label{eq-kernel2}
& = & 
\ds\frac{A_n^{\delta-1}}{A_n^\delta}\int_{-1}^1k_n^{\delta-1}(z(\bdx;t))\, 
w^{(\alpha+\beta+\frac{1}{2},\beta)}(t)\, dt,
\end{eqnarray}
where $k_n^\delta$ is the univariate Ces\`aro $(C,\delta)$ kernel, 
$$
k_n^\delta(t) = 
\sum_{k=0}^n \frac{A_{n-k}^\delta}{A_n^\delta}\,
h_k^{(\alpha+\beta+\frac 12,\beta)}
P_k^{(\alpha+\beta+\frac 12,\beta)}(1)
P_k^{(\alpha+\beta+\frac 12,\beta)}(t).
$$
\end{cor}

\medskip
Equation (\ref{eq-kernel2}) allows us to establish (\ref{L1norm}) by working
with the univariate kernel $k_n^\delta$. This leads to our main result in this
paper.

\begin{theorem}\label{th-main}
Let $\alpha-\frac 12\geq\beta\geq 0$. The $(C,\delta)$ means of the 
orthogonal expansions with respect to $W_B^{\alpha,\beta}$ converge 
uniformly to $f$ on $B$ for every continuous function $f$ if  
$\delta>\alpha+\beta+1$. 
\end{theorem}

The proof of Theorem \ref{th-main} will be stated in Section \ref{sec:proof}. 
The proof involves sharp estimates on various pieces, which indicates that
the order $\delta >\alpha+\beta+1$ is sharp. In other words, we conjecture that
the convergence fails if $\delta \le \alpha + \beta +1$. 

The compact formula (\ref{eq-kernel2}) also allows us to deduce
positivity of the $(C,\delta)$ means on $B$ from its counterpart of
univariate Jacobi expansions.  

\begin{theorem}\label{th-positivity}
Let $\alpha-\frac 12\ge\beta\ge 0$. The $(C,\alpha+2\beta + 3/2)$ means of
the orthogonal expansions with respect to $W_B^{\alpha,\beta}$ define  
positive linear operators. 
\end{theorem}

\medskip
{\it Proof.} 
In \cite{Ga77} it is proved that the $(C,\mu+\nu+2)$ means of the 
univariate Jacobi expansion 
$$
\sum_{n=0}^\infty [h_n^{\mu,\nu}]^{-1} P_n^{(\mu,\nu)}(1)
P_n^{(\mu,\nu)}(x)
$$ 
are nonnegative for $-1\le x\le 1$, $\mu\ge\nu \ge -1/2$. 
By Corollary \ref{kernel-reduced}, the $(C,\delta)$ means of 
$K_n(\bdx)$ are the integrals of the $(C,\delta-1)$ means of the Jacobi 
expansions for $w^{(\alpha+\beta+\frac12,\beta)}$, which shows that
$K_n^\delta(\bdx)$ is nonnegative for $\delta = \alpha+2\beta + 3/2$. 
By the product formula, $\mathcal{K}_n^\delta(\bdx,\bdy) = T_{\bdx}
 \mathcal{K}_n(\mathbf{e}, \mathbf{y}),$ and $T_{\bdx}$ is an integral
 against a nonnegative measure, it follows that 
$\mathcal{K}_n^\delta(\bdx,\bdy) 
 \ge 0$ for $\delta  = \alpha+2\beta + 3/2$. 
\hfill{$\Box$}

\medskip
We note that if the $(C,\delta_0)$ means are nonnegative, then the $(C,\delta)$ 
means are nonnegative for all $\delta\geq\delta_0$.  The positivity implies 
$\|K_n^\delta\|_1 = 1$, hence, convergence of the means. However, since 
convergence of the $(C,\delta)$ means implies the convergence of the $(C,\delta')$ 
means for $\delta' > \delta$, the  convergence also follows form 
Theorem 3.2.

\section{Proof of Theorem \ref{th-main}}
\label{sec:proof}

We start with a  result in \cite[p. 261, (9.4.13)]{Sz75} and its extension in 
\cite{Li95} given in the following lemma.

\begin{lem} 
\label{lem:3.6}
Let $k_n^\delta(w^{(\xi,\eta)};u)$, $u\in[-1,1]$, denote the
kernel for the univariate Ces\`aro $(C,\delta)$ means of the
Jacobi expansion with respect to the weight function $w^{(\xi,\eta)}$
on $[-1,1]$. Then for any $\xi,\eta > -1$ such that 
$\xi+ \eta+ \delta + 3 >0$ we have that
$$
  k_n^\delta(w^{(\xi,\eta)};t) = \sum_{j=0}^J b_j(\xi,\eta,\delta,n)
    P_n^{(\xi+\delta+j+1,\eta)}(t) + G_n^{\delta}(t),
$$
where $J$ is a fixed integer and
$$
   G_n^\delta(t) = \sum_{j=J+1}^\infty d_j(\xi,\eta,\delta,n)
  k_n^{\delta+j}(w^{(\xi,\eta)}, 1,t).
$$
Moreover, the coefficients satisfy the inequalities
$$
 |b_j(\xi,\eta,\delta,n)| \le c n^{\xi+1-\delta - j} \quad \hbox{and}
 \quad  |d_j(\xi,\eta,\delta,n)| \le c j^{-\xi-\eta-\delta - 4}.
$$
\end{lem}

Furthermore, we will need an estimate for the univariate kernel
which was proved in \cite[Lemma 5.2]{Li03}.
\begin{lem} 
\label{lem:5.2}
Let $k_n^\delta(w^{(\xi,\eta)};u)$, $u\in[-1,1]$, denote the
kernel for the univariate Ces\`aro $(C,\delta)$ means of the
Jacobi expansion with respect to the weight function $w^{(\xi,\eta)}$
on $[-1,1]$. Let further $\xi,\eta > -1$ and
$\delta\geq \xi+\eta+2$. Then
$$
|k_n^\delta(w^{(\xi,\eta)};t)| \, \leq \,
cn^{-1}(1-t-n^{-2})^{-(\xi+3/2)}.
$$
\end{lem}

It is well known that the Jacobi polynomials satisfy the following
estimate (\cite[(7.32.5) and (4.1.3)]{Sz75}).
\begin{lem} \label{lem:3.2}
For $\alpha,\beta>-1$ and $t \in [0,1]$,
\begin{equation} 
|P_n^{(\alpha,\beta)} (t)| \le c n^{-1/2} (1-t+n^{-2})^{-(\alpha+1/2)/2}.
\end{equation}
The estimate on $[-1,0]$ follows from the fact that $P_n^{(\alpha,\beta)}
(t) = (-1)^nP_n^{(\beta, \alpha)} (-t)$.
In particular, for all $t \in [-1,1]$, we have the estimate
\begin{equation} \label{Est-Jacobi}
|P_n^{(\alpha,\beta)} (t)| \le c n^{-1/2} 
\left[(1-t+n^{-2})^{-(\alpha+1/2)/2}
     + (1-t+n^{-2})^{-(\b+1/2)/2}
\right].
\end{equation}
\end{lem}

The central part of the proof is the following proposition. 
\begin{pro}
\label{proposition}
If $\d > \a + \b +1$ then 
$$
 \int_B \int_{-1}^1 |P_n^{(\a+\b+\d+\frac12,\b)}(z(\bdx;t))| 
   w^{(\a+\b+\frac12,\b)}(t)\, dt \,
   W^{\a,\b}(\bdx)\, d\bdx \le c n^{\d - \a- \b-3/2}.  
$$
\end{pro}

{\it Proof. }
Using the inequality (\ref{Est-Jacobi}), we see that it is sufficient 
to show that 
\begin{equation} 
\label{Est1}
  J_1:= \int_B \int_{-1}^1 \frac{w^{(\a+\b+\frac12,\b)}(t)}
     {(1- z(\bdx;t) + n^{-2})^{\frac{\a+\b+\d+1}{2}}}\, dt \, 
         W^{\a,\b}(\bdx)\, d\bdx \le c n^{\d - \a- \b-1},
\end{equation}
and 
\begin{equation} 
\label{Est2}
  J_2:= \int_B \int_{-1}^1 \frac{w^{(\a+\b+\frac12,\b)}(t)}
     {(1+ z(\bdx;t) + n^{-2})^{\frac{\b+1/2}{2}}}\, dt \, 
         W^{\a,\b}(\bdx)\, d\bdx \le c n^{\d - \a- \b -1}.
\end{equation}

We start with $J_1$. Its estimate is divided into several cases, 
according to the decompositions $[-1,1]=[-1,0]\cup[0,1]$ and 
$$
B=B_+\cup B_-, \quad  B_+ :=\{\bdx=(x_1,x_2) \in B: x_1\geq 0\}, \quad 
         B_-:=\{\bdx=(x_1,x_2) \in B: x_1\leq 0\}.
$$ 
To simplify notation, we further denote 
$$
    \gamma : = (\a+\b+\d+1)/2
$$
throughout this proof.  The following basic identity can be easily 
verified,
\begin{equation} 
\label{1-z}
  1- z(\bdx;t) =  (1-t^2)(1-x_1) + \frac12 (1-t)^2 (1-x_2). 
\end{equation}

\medskip\noindent
{\it Case 1.} The integral over $\bdx\in B$ and $t\in [0,1]$.

Since $t \ge 0$, by (\ref{1-z}), $1-z(\bdx;t) \ge (1-t)(1-x_1)$, and  
$w^{(\a+\b+\frac12,\b)}(t)\le c (1-t)^{\a + \b +1}$, so that 
$$
 \int_0^1 \frac{w^{(\a+\b+\frac12,\b)}(t)}
     {(1- z(\bdx;t) + n^{-2})^\g} dt 
  \le  \int_0^1 \frac{ (1-t)^{\a + \b +1} }
     {[(1- t)(1-x_1) + n^{-2}]^\g} dt =: f(x_1).
$$
Since $f$ depends on $x_1$ only, it readily follows that 
$$
\int_{B}f(x_1)W^{\alpha,\beta}(\bdx)\,d\bdx 
 =  \int_{-1}^1  f(x_1) \int_{x_1^2}^1W^{\alpha,\beta}(\bdx) dx_2 dx_1
 =  c \int_{-1}^1f(x_1)(1-x_1^2)^{\alpha+\beta} dx_1, 
$$
where we used the fact that
$$
\int_{x_1^2}^1(1-x_2)^{\alpha-1/2}(x_2-x_1^2)^{\beta-1/2}\ dx_2
\, = \,
\frac{\Gamma(\alpha+\frac{1}{2})\Gamma(\beta+\frac{1}{2})}
     {\Gamma(\alpha+\beta+1)}(1-x_1^2)^{\alpha+\beta}.
$$
Thus, it follows that this part of the integral in $J_1$ is bounded by
$$
  \int_B \int_0^1 \,\cdots\,
  \le c \int_{-1}^1 \int_0^1   \frac{(1-t)^{\a+\b+\frac12} }
     {[(1- t)(1-x_1) + n^{-2}]^\g} dt (1-x_1^2)^{\alpha+\beta}  dx_1
   : = I_1^+ + I_1^-,
$$
where we split of the integral over $[-1,1]$ as two integrals over $[0,1]$
and $[-1,0]$, respectively, and define $I_1^+$ and $I_1^-$ accordingly. 
For $I_1^+$,  we have
\begin{eqnarray*}
I_1^+ 
& = &
\int_0^1\int_0^1
\frac{s^{\alpha+\beta+\frac{1}{2}}(1-x^2)^{\alpha+\beta}}
{[s (1-x)+n^{-2}]^\g} ds \,dx 
  = \int_0^1\int_0^x\frac{u^{\alpha+\beta+\frac{1}{2}}
(2+x)^{\alpha+\beta}}{(u+n^{-2})^\g} du \, x^{-\frac{3}{2}} dx \\
& \leq &
  \int_0^1   \frac{u^{\alpha+\beta+\frac{1}{2}} } {(u+n^{-2})^\g}
     \int_u^1 x^{-\frac32} dx\,  du 
  \le c \int_0^1  \frac{u^{\alpha+\beta} } {(u+n^{-2})^\g} du \\
& = &
c n^{ 2 \g - 2 \a - 2 \b -2} \int_0^{n^2} 
\frac{v^{\alpha+\beta} } {(1+n^{2}v)^\g} dv 
   \le c n^{\d-\a-\b-1},
\end{eqnarray*}
as the last integral is bounded if $\g > \a + \b + 1$.  
For $I_1^-$, we have $1- x_1 \ge 1$, so that 
$$
  I_1^- \le c\int_0^1   \frac{(1-t)^{\a+\b+\frac12} }
     {[(1- t)  + n^{-2}]^\g} dt  \int_{-1}^0 (1-x_1^2)^{\alpha+\beta}  dx_1
  \le c  \int_0^1  \frac{u^{\alpha+\beta+1} } {(u+n^{-2})^\g} du, 
$$
which has the same bound as $I_1^+$ if we use $u^{\a+\b+1} \le u^{\a+\b}$. 

\medskip\noindent
{\it Case 2.} The integral over $\bdx\in B_+$ and $t\in [-1,0]$. 

For $t \le 0$, by (\ref{1-z}), $1-z(\bdx;t) \ge (1+t)(1-x_1) + \frac12 (1-x_2)$ 
and $w^{(\a+\b+1,\b)}(t) \le c (1+t)^\b$.  Hence, this portion of the integral 
in $J_1$ is bounded by
$$
 I: = c \int_{B_+} \int_{-1}^0  \frac{(1+t)^\beta} 
    {\big[(1+t)(1-x_1)+(1-x_2)+n^{-2}\big]^{\g}}\,dt\, W^{\alpha,\beta}(\bdx)\,d\bdx.
$$
Changing variables $u = (1+t)(1-x_1)$ shows that 
\begin{eqnarray*}
I 
& = &
   c \int_{B_+} \int_{0}^{1-x_1} \frac{ u ^\beta} 
    {\big[ u +(1-x_2)+n^{-2}\big]^{\g}}\,du\, 
   \frac{W^{\alpha,\beta}(\bdx)}{(1-x_1)^\b}\,d\bdx \\
& = &
   c \int_0^1 \int_0^{\sqrt{x_2}} \int_0^{x_1}\frac{ u ^\beta} 
    {\big[ u +(1-x_2)+n^{-2}\big]^{\g}}\,du 
      \frac{(x_2-x_1^2)^{\beta-\frac 12}} {(1-x_1)^{\beta+1}}dx_1 
          (1-x_2)^{\alpha-\frac 12} dx_2. 
\end{eqnarray*}
Let us first consider the two inner integrals. Changing the order of 
integration shows that 
$$
\int_0^{\sqrt{x_2}}\int_0^{1-x_1} du dx_1 =
\int_0^{1-\sqrt{x_2}} \int_0^{\sqrt{x_2}} d x_1 du +
 \int_{1-\sqrt{x_2}}^1 \int_0^{1-u} d x_1 du.
$$
Now, using $x_2 + \sqrt{x_1} \sim \sqrt{x_2}$ and $\beta+\frac12 > 0$, 
integration by parts once gives, 
\begin{eqnarray*}
&& \int_0^{\sqrt{x_2}}  \frac{(x_2-x_1^2)^{\beta-\frac 12}}
     {(1-x_1)^{\beta+1}}dx_1 
 \le c  x_2^{\b -\frac12} \int_0^{\sqrt{x_2}} 
     \frac{(\sqrt{x_2}-x_1)^{\beta-\frac 12}}
       {(1-x_1)^{\beta+1}}dx_1 \\
&& \, = c  \sqrt{x_2}^{\b -\frac12} 
     \left[\frac{\sqrt{x_2}^{\b + \frac12} }{\b+\frac12} 
         + \frac{\b+1}{\b+\frac12} \int_0^{\sqrt{x_2}}  
           \frac{1} {(1-x_1)^{\frac 32}}\,dx_1 
     \right]  
    \le  c \frac {\sqrt{x_2}^{\b -\frac12}}{(1-x_2)^\frac 12}.
\end{eqnarray*}
Analogously, using $\sqrt{x_2} \le \sqrt{x_2}+ x_1 \le 2$ and 
integration by parts once, we have the estimate
\begin{eqnarray*}
\int_0^{1-u} \frac{(x_2-x_1^2)^{\beta-\frac 12}}{(1-x_1)^{\beta+1}}\,dx_1
& \leq &
c \max\{1, \sqrt{x_2}^{\b -\frac12} \}
       \int_0^{1-u} \frac{(\sqrt{x_2}-x_1)^{\beta-\frac 12}}
          {(1-x_1)^{\beta+1}}\,dx_1 \\
& \leq &
c \max\{1, \sqrt{x_2}^{\b -\frac12} \}
         \frac{1}{\sqrt u}.
\end{eqnarray*}
Adding the two parts together, we obtain that for a generic 
function $f(u)$, 
\begin{eqnarray*}
&& \int_0^{\sqrt{x_2}}\int_0^{1-x_1} f(u) du
\frac{(x_2-x_1^2)^{\beta-\frac 12}}
     {(1-x_1)^{\beta+1}} dx_1 \\
&& \qquad \le c \max\{1, \sqrt{x_2}^{\b -\frac12} \}
    \int_0^1 f(u) \max\{(1-x_2)^{-\frac 12},u^{-\frac 12}\} du.
\end{eqnarray*} 
Consequently, we conclude that 
$$
  I  \le c  \int_0^1  \int_0^1 
   \frac{ u ^\beta  \max\{1, \sqrt{x_2}^{\b -\frac12} \} } 
        {\big[ u +(1-x_2)+n^{-2}\big]^{\g}}du 
      \left[ (1-x_2)^{-\frac12} + u^{\frac12} \right] 
      (1-x_2)^{\alpha-\frac 12} dx_2. 
$$
We note that the term $\max\{1,\sqrt{x_2}^{\b - \frac12}\}$, which 
matters only if $\b < \frac12$,  is integrable as $\b > -1/2$ and 
it plays a minor role. In fact, if we split the integral of $x_2$ as 
an integral over $[0, \frac12]$ and $[\frac12, 1]$, then the part 
over $x_2 \in [0,\frac12]$ is bounded as 
$u + (1-x_2) + n^{-2} \ge \frac12$. Hence, we only need to estimate
the sum
$$
  \int_0^1\int_0^1
\frac{u^\beta x^{\alpha-1}}{(u+x+n^{-2})^\gamma}du\,dx + 
  \int_0^1\int_0^1
\frac{u^{\beta-\frac 12} x^{\alpha-\frac 12}}{(u+x+n^{-2})^\gamma}du\,dx 
:= I_1 +I_2.
$$
For $I_1$ we change variables $v= u+x$ and then exchange the order of 
the integrals, 
\begin{eqnarray*}
I_1  = \int_0^1\int_x^{x+1}
   \frac{(v-x)^\beta x^{\alpha-1}}{(v+n^{-2})^\gamma}\,dv\,dx    
& = &
  \int_0^1 \frac 1{(v+n^{-2})^\gamma}
  \int_0^v (v-x)^\beta x^{\alpha-1}\,dx\,dv  \\
& + &  
  \int_1^2 \frac 1{(v+n^{-2})^\gamma}
  \int_{v-1}^1 (v-x)^\beta x^{\alpha-1}\,dx\,dv.
\end{eqnarray*}
The second term is bounded by a constant as $v \ge 1$, whereas the inner 
integral of the first term is a Beta integral and equals 
$v^{\beta+\alpha} B(\beta+1,\alpha)$, so that 
$$
I_1 \le c 
    \left( 1+\int_0^1 \frac{v^{\alpha+\beta}}{(v+n^{-2})^\gamma}\,dv 
    \right)
    \le c n^{\d - \a -\b -1},
$$
where the second inequality follows from the last step in the estimate for 
$I_1^+$ of Case 1. Notice that $I_2$ becomes $I_1$ if we replace $(\a,\b)$ 
by $(\a-\frac12, \beta+\frac12)$, so that $I_2$ has the same upper bound.

\medskip\noindent
{\it Case 3.} The integral over $\bdx\in B_-$ and $t\in [-1,0]$. 

Just like in Case 2, this portion of the integral in $J_1$ is bounded by
$$
 I: = c \int_{B_-} \int_{-1}^0  \frac{(1+t)^\beta} 
    {\big[(1+t)(1-x_1)+(1-x_2)+n^{-2}\big]^{\g}}\, dt\, 
    W^{\alpha,\beta}(\bdx)\,d\bdx.
$$
For $\bdx\in B_-$, we have that $1-x_1 \ge 1$ so that we can drop 
the factor $1-x_1$ in the denominator. Changing variables $u = 1+t$ then 
shows that
\begin{eqnarray*}
  I 
& \le &
 c  \int_{B_-} \int_{-1}^0  \frac{u^\beta} 
       {\big[ u +(1-x_2)+n^{-2}\big]^{\g}}\,du 
       (1-x_2)^{\a - \frac12}
       (x_2 - x_1^2)^{\b-\frac12}\,d\bdx \\
& \le &
 c \int_0^1\int_0^1 \frac{u^\beta\, du}
     {(u+(1-x_2)+n^{-2})^\gamma} (1-x_2)^{\alpha-\frac 12} 
   \int_{-\sqrt{x_2}}^0 (x_2-x_1^2)^{\beta-\frac 12}\,dx_1\,dx_2\\
& \le &
 c  \int_0^1\int_0^1  \frac{u^\beta  (1-x_2)^{\alpha -1/2}}
     {(u+(1-x_2)+n^{-2})^\gamma}  x_2^{\b} \,du\,  dx_2
      \le c  \int_0^1\int_0^1  \frac{u^\beta  x^{\alpha -1/2}}
     {(u+x+n^{-2})^\gamma} \,du\,  dx.  
\end{eqnarray*}
Comparing to $I_2$ in Case 2 and using $u^\b \le u^{\b-\frac12}$, 
we see that $I$ is bounded by $c n^{\d-\a-\b-1}$ as before. 

Putting these cases together completes the proof of (\ref{Est1}).

We now prove (\ref{Est2}).  A straightforward computation shows that 
$$
1+ z(\bdx;t)  =  
\frac 12 \left[ \big(1+t+(1-t)x_1\big)^2 +(1-t)^2(x_2-x_1^2) \right] \ge 
\frac 12 (1-t)^2(x_2-x_1^2). 
$$
Using this inequality in $J_2$, we obtain that
\begin{eqnarray*}
J_2 
& \le &
 c \int_B \int_{-1}^1
  \frac{(1-t)^{\alpha+\beta+\frac 12}(1+t)^\beta}
     {[(1-t)\sqrt{x_2-x_1^2}+n^{-1}]^{\beta+\frac 12}}\,dt 
      (1-x_2)^{\alpha-\frac12}(x_2 - x_1^2)^{\beta-\frac12}\, d\bdx \\
& \le &
 c \int_B\int_{-1}^1
  \frac{\big((1-t)\sqrt{x_2-x_1^2}\big)^{\beta+\frac 12}}
     {[(1-t)\sqrt{x_2-x_1^2}+n^{-1}]^{\beta+\frac 12}}\,
      (1-t)^{\alpha+\frac 12}(1+t)^\beta \, dt \\
& &
 \qquad\qquad \qquad\qquad \cdot
(1-x_2^2)^{\alpha-\frac 12}(x_2-x_1^2)^{\frac{\beta-\frac 32}2}\,
 d\bdx \\
& \le &
 c \int_B\int_{-1}^1 (1-t)^{\alpha+\frac 12}(1+t)^\beta\, dt
(1-x_2^2)^{\alpha-\frac 12}(x_2-x_1^2)^{\frac{\beta-\frac 32}2}\,
 d\bdx \le c. 
\end{eqnarray*}
Hence, for $\d > \a + \b +1$, $J_2$ is bounded and this completes the
proof. 
\proofend

{\it Proof of Theorem \ref{th-main}. }
By Corollary \ref{cor2.2}, it is sufficient to show that 
$$
\|K_n^\delta\|_1= \int_B |K_n^\delta(\bdx)| W^{\alpha,\beta}(\bdx)\, d\bdx
\, \leq \, c
$$
under the condition that $\delta > \alpha+ \beta+1$. We set 
$J=\alpha+2\beta+\frac 32$ in Lemma \ref{lem:3.6} and obtain that
\begin{eqnarray*}
|k_n^\delta(t)| 
& \leq &
\sum_{j=0}^J |b_j(\alpha+\beta+\tfrac{1}{2},\beta,\delta,n)
P_n^{(\alpha+\beta+\delta+j+\frac 32,\beta)}(t)| \\
& & + \,
\sum_{j=J+1}^\infty |d_j(\alpha+\beta+\tfrac{1}{2},\beta,\delta,n)
k_n^{\delta+j}(t)|.
\end{eqnarray*}
Using Corollary \ref{kernel-reduced} together with Proposition \ref{proposition}, 
and  taking into account the estimate for the coefficients given in 
Lemma \ref{lem:3.6}, we obtain that
\begin{eqnarray*}
\|K_n^\delta\|_1\,
& \leq &
\frac{A_n^{\delta-1}}{A_n^{\delta}}\Big[
cn\, + \,
\sum_{j=J+1}^\infty |d_j(\alpha+\beta+\tfrac{1}{2},\beta,\delta,n)| \\
& & \times \,
\int_B\int_{-1}^1 
|k_n^{\delta+j}(z(\bdx;t))| w^{(\alpha+\beta+\frac12,\beta)}(t) 
\,dt\, W^{\alpha,\beta}(\bdx)\, d\bdx
\Big].
\end{eqnarray*}
To estimate the second sum, we use Lemma \ref{lem:5.2}. Thus,
we have to derive an upper bound for the integral
$$
I:=
n^{-1}\,
\int_B\int_{-1}^1 
\frac{w^{(\alpha+\beta+\frac 12,\beta)}(t)}
     {(1-z(\bdx;t)+n^{-2})^{\alpha+\beta+2}}
dt\, W^{\alpha,\beta}(\bdx)\,d\bdx.
$$
Setting $\gamma=\alpha+\beta+2$, we already proved 
(cf. (\ref{Est1})) that
$$
\int_B\int_{-1}^1
\frac{w^{(\alpha+\beta+\frac 12,\beta)}(t)}
     {(1-z(\bdx;t)+n^{-2})^\gamma}\, dt\,
W^{\alpha,\beta}(\bdx)\, d\bdx
\, \leq \,
cn^{2\gamma-2\alpha-2\beta-2} = c n^2,
$$
provided that $\gamma>\alpha+\beta+1$. Hence, we obtain the bound
$I \le c n$. From 
${{n+\delta-1}\choose n}/{{n+\delta}\choose n}=\frac \delta{n+\delta}$
it follows that $\|K_n^\delta||_1$ is bounded uniformly in $n$.
\proofend


\small
Wolfgang zu Castell \& Frank Filbir\\
Institute of Biomathematics and Biometry\\
Helmholtz Zentrum M\"{u}nchen\\
German Research Center for Environmental Health\\
Ingolstaedter Landstrasse 1, 85764 Neuherberg, Germany\\
{\tt [castell,filbir]@gsf.de}

\smallskip
Yuan Xu\\
Department of Mathematics
University of Oregon,\\
Eugene, OR 97403-1222, U.S.A.\\
{\tt yuan@math.uoregon.edu}

\end{document}